\documentclass[10pt,a4paper,twoside]{article}

\usepackage{amssymb,amsmath,multirow}

\newtheorem{theorem}{Theorem}[section]

\newtheorem{example}[theorem]{Example}

\newtheorem{lemma}[theorem]{Lemma}

\newcommand{\integ}{\int\limits_1^\infty}

    \makeatletter
    \let\@fnsymbol\@arabic
    \makeatother

\begin{document}

\title{Goodness of Fit Tests for Pareto Distribution Based on a Characterization and their Asymptotics}

\author{  Marko Obradovi\' c \footnote{marcone@matf.bg.ac.rs} , Milan Jovanovi\'c\footnote{mjovanovic@matf.bg.ac.rs} ,
        Bojana Milo\v sevi\'c\footnote{bojana@matf.bg.ac.rs} \\\medskip
{\small Faculty of Mathematics, University of Belgrade, Studenski trg 16, Belgrade, Serbia}}

\date{}
\maketitle

\begin{abstract}
In this paper we present a new characterization of Pareto distribution and  consider  goodness of fit tests based on it.
We provide an integral  and Kolmogorov- Smirnov type statistics based on U-statistics and we calculate Bahadur efficiency for various alternatives.
We find locally optimal alternatives for those tests. For small sample sizes we compare the power of those tests
with some common goodness of fit  tests.

\end{abstract}

{\small \textbf{ keywords:} Pareto distribution, Bahadur efficiency, U-statistics, Goodness of fit test

\textbf{MSC(2010):} 62G10, 62G20, 62H05}

\section{Introduction}
A very important part of data analysis is ensuring that the data come from a particular family of distributions.
In order to achieve this goal, there exist a variety of goodness of fit tests.
Some of such tests are based on empirical distribution function.

A characterization of a certain family of distributions is a property which is true only for that family.
See \cite{Galambos} for more on characterizations. The characterizations are often
the best way to distinguish one family of distributions from the others.
Hence, they are convenient for use in goodness of fit testing. Such tests are often free of some parameters.
Some examples of such tests can be found in \cite{Angus}, \cite{Morris}, \cite{Noughabi}.

One of the basic problems in hypotheses testing is choosing the more efficient test of two proposed. 
One way to compare them is the asymptotic relative efficiency. Several approaches have been used. 
The Bahadur efficiency has the advantage that it can also be applied to the statistics with non-normal 
asymptotic distribution such as Kolmogorov-Smirnov, Cramer-von Mises and others.

In particular, the Bahadur asymptotic efficiency  of tests based on characterization have been discussed in numerous papers, see e.g.
\cite{Nikitin4}.
The  Bahadur asymptotic efficiency  of a test for exponentiality based on characterization has been studied in \cite{Nikitin3}.
Integral and Kolmogorov-Smirnov type test statistics for testing exponentiality based on characterization and their  Bahadur asymptotic efficiency
have been studied in \cite{Volkova}.
Tests of the same kind for power function distribution have been studied in \cite{Nikitin2}.

The family of Pareto distributions is frequently used in many fields to model heavy tailed phenomena. 
Initially, its application was in the modelling income over a population, see \cite{Asimit}.
It has been used for modelling socio-economic quantities such as the size of cities and 
the firm sizes, as in \cite{Simon1}, \cite{Simon2}. In finance, it found its application in modelling stock price returns, see
\cite{Rydberg},
while in actuaries it is a frequent distribution of quantities such as
 excess of loss quotations in insurance, as in \cite{Rytgaard}. In recent times, it was found suitable for internet teletraffic modelling,
 e.g. file sizes and FTP transfer times, as in \cite{Chlebus}.
 There are also many of its applications in physics, hydrology and seismology, see e.g. 
 \cite{Aban}.

In most of these applications the quantities in question are generated by stochastic processes with limiting distribution
which is either the Pareto distribution or some other distribution strongly
skewed to the right. Therefore, it is relevant to test the hypothesis that data come from the Pareto distribution against 
such alternatives.

Some goodness of fit tests for Pareto distribution can be found in \cite{Gulati}, \cite{Martynov}, \cite{Rizzo}. 
Some characterizations of Pareto distributions can be found in \cite{Hitha}.

In this paper we propose two goodness of fit tests based on a new characterization of Pareto distribution.
In section 2 we introduce the characterization. We propose two test statistics, $T_n$ of integral and $V_n$ of Kolmogorov-Smirnov type.
We study the statistic $T_n$ in section 3. Using U-statistic we investigate asymptotic behaviour of $T_n$. We calculate its local
Bahadur efficiency with respect to several alternatives. We also find some of locally optimal alternatives. Using Monte-Carlo methods we
calculate critical values of this test.
In section 4 we do the analogous study for the statistic $V_n$.
In section 5 we compare the powers of our tests with some standard goodness of fit tests.

\section{The Test Statistics Based on New Characterization of Pareto Distribution}

Let $\mathcal{P}$ be the family of Pareto distributions  with the distribution function
\begin{equation*}
F(x)=1-x^{-\alpha}, x\geq 1, \;\; \alpha>0.
\end{equation*}

We now present new characterization.
\begin{theorem}
 Let $X$ and $Y$ be i.i.d. non-negative absolutely continuous  random variables. Then,  $X$ and $\max\{\frac{X}{Y},\frac{Y}{X}\}$
have the same distribution if and only if the distribution of $X$ belongs to the family $\mathcal{P}$.
 \end{theorem}

 \noindent\textbf{Proof.} Let $\widetilde{X}=\ln X$ and $\widetilde{Y}= \ln Y$. Then

 \begin{equation*}
  \max\{\frac{X}{Y},\frac{Y}{X}\}=\max\{e^{\widetilde{X}-\widetilde{Y}},e^{\widetilde{Y}-\widetilde{X}}\}=e^{|\widetilde{X}-\widetilde{Y}|}
 \end{equation*}
 and
 \begin{equation*}
  \ln(\max\{\frac{X}{Y},\frac{Y}{X}\})=|\widetilde{X}-\widetilde{Y}|.
 \end{equation*}

 Since the logarithm is a monotonous transformation, then the statement that $X$ and $\max\{\frac{X}{Y},\frac{Y}{X}\}$
 have the same distribution  is equivalent to the statement
 that $\widetilde{X}$ and  $|\widetilde{X}-\widetilde{Y}|$ are identically distributed.
 In \cite{PuriRubin} it was proven that the only continuous distributions
 that satisfy this are exponential distributions, so $\widetilde{X}$
 must have exponential distribution with some scale parameter $\alpha$. Since $X=e^{\widetilde{X}}$, then $X$ has Pareto
 distribution with the same parameter $\alpha$. $\hfill \Box$
 \medskip

 Let $(X_1,X_2,\ldots, X_{n})$ be a sample from non-negative continuous distribution function $F$.
 Using above characterization, we are going to construct two goodness of fit tests for Pareto distribution.
 Since the characterization
 is free of the shape parameter $\alpha$, we shall test the composite null hypothesis $H_0:F \in \mathcal{P}$, against the general alternative
 $H_1:F\notin \mathcal{P}$ with the same support $[1,\infty)$.

  The test statistics we are going to use  are
 \begin{equation*}
       T_n=\int\limits_{1}^\infty(M_n(t)-F_n(t))dF_n(t)
   \end{equation*}
   and

 \begin{equation*}
  V_n=\sup\limits_{t\geq 1}|M_n(t)-F_n(t)|,
 \end{equation*}
where $F_n(t)=\frac{1}{n}\sum\limits_{i=1}^{n}{\rm I}\{X_{i}\leq t\}$ is the empirical distribution function based on the sample
 $(X_1,X_2,\ldots, X_n)$ and $M_n(t)$ is  $U$-empirical distribution function based on above characterization defined by
 \begin{equation*}
  M_n(t)=\binom{n}{2}^{-1}
  \sum\limits_{i=1}^{n-1}\sum\limits_{j=i+1}^n{\rm I} \left\{\max\left(\frac{X_i}{X_j},\frac{X_j}{X_i}\right)\leq t\right\},\;\;t\geq 1.
 \end{equation*}
 
 We shall assume without loss of generality that large values of test statistics are significant.

 \section{Test Statistic $T_n$}

  In this section we examine the properties of integral type statistic $T_n$. Using asymptotical equivalence of this statistic
  and the $U$-statistic with kernel
\begin{eqnarray*}
 \Upsilon(X,Y,Z)&=&\frac{1}{3}\;{\rm I}\left\{\max\left(\frac{X}{Y},\frac{Y}{X} \right)\leq Z\right\}
 +\frac{1}{3}\;{\rm I}\left\{\max\left(\frac{X}{Z},\frac{Z}{X} \right)\leq Y\right\}\\&+&
 \frac{1}{3}\;{\rm I}\left\{\max\left(\frac{Y}{Z},\frac{Z}{Y} \right)\leq X\right\}-\frac{1}{2},
\end{eqnarray*}
by the law of large numbers for $U$-statistics (see \cite{serfling})
we get that
\begin{equation}\label{tnconv}
 T_n\overset{p}{\underset{n\rightarrow\infty}{\longrightarrow}}P\left\{\max\left(\frac{X}{Y},\frac{Y}{X} 
  \right)\leq Z\right\}-\frac{1}{2}.
\end{equation}

Let us now examine the asymptotic behaviour of $T_n$ under null hypothesis via $U$-statistic defined above.

Let $\upsilon(X)$ be the projection of $\Upsilon(X,Y,Z)$ on $X$. Then
\begin{eqnarray*}
\nonumber \upsilon(s)&=&E(\Upsilon(X,Y,Z)|X=s)\\&=&
\frac{2}{3}P\left\{\max\left(\frac{s}{Y},\frac{Y}{s}\right)\leq Z \right\}
+\frac{1}{3}P\left\{\max\left(\frac{Y}{Z},\frac{Z}{Y}\right)\leq s \right\}-\frac{1}{2}\\
 &=&\frac{2}{3}(\int\limits_{1}^sdz\int\limits_{\frac{s}{z}}^\infty
 \alpha^2z^{-\alpha-1}y^{-\alpha-1}dy+\int\limits_{s}^\infty dz\int\limits_{\frac{z}{s}}^\infty\alpha^2z^{-\alpha-1}y^{-\alpha-1}dy)\\
&+& \frac{1}{3}P\left\{Y\leq s\right\}-\frac{1}{2}
\\&=&\frac{2}{3} \frac{2\alpha\ln s+1}{2s^\alpha}+\frac{1}{3}(1-s^{-\alpha})-\frac{1}{2}
\\&=&\frac{2}{3}\frac{\alpha \ln s}{s^{\alpha}}-\frac{1}{6}.
\end{eqnarray*}
It is easy to see that $E(\upsilon(X))=0$, so
the variance of this projection is
\begin{equation}\label{sigma}
 \sigma^2={\rm Var}(\upsilon(X))=\int\limits_{1}^\infty \upsilon^2(s)\alpha s^{-\alpha-1}ds=\frac{5}{972}.
\end{equation}
Since  the variance of the projection is positive, the kernel $\Upsilon(X,Y,Z)$ is not degenerate,
so we can apply Hoeffding's theorem for $U$-statistics with non-degenerate kernels  (see \cite{Hoeffding}).
Since the degree  of $U$-statistics is 3, the asymptotic variance is $3^2\sigma^2$, so we obtain
\begin{equation}\label{asydis}
 \sqrt{n}T_n\overset{d}{\rightarrow}\mathcal{N}\left(0,\frac{5}{108}\right).
\end{equation}

\subsection{Bahadur efficiency}

This way of measuring asymptotic efficiency is explained in detail in \cite{Bahadur}, \cite{Nikitin1}. For two tests with the
same null ($H_0:\theta\in\Theta_0$) and alternative ($H_1:\theta\in\Theta_1$) hypotheses, the asymptotic
relative Bahadur efficiency is defined as the ratio of sample sizes needed to reach the same test power when the level of significance
approaches zero.
It can be expressed as the ratio of
Bahadur exact slopes, functions proportional to exponential rate for a sequence of test statistics, provided that these 
functions exist. This is fulfilled for the majority of tests.

The Bahadur exact slope (see \cite{Nikitin1}) can be evaluated as

\begin{equation}\label{slope}
 c_{T}(\theta)=2f(b_{T}(\theta)),
\end{equation}
where $T_{n}\overset{p}{\underset{n\rightarrow \infty}{\longrightarrow}} b_{T}(\theta)$ for $\theta \in \Theta_1$,
$G_{n}(t)=\inf\{P_{\theta}\{T_n<t\},\theta\in \Theta_0\}$ and
$f(t)=-\lim\limits_{n\rightarrow \infty} \frac{1}{n}\ln(1-G_{n}(t))$ for each $t$ from an open interval $I$ on which $f$ is continuous
and $\{b_T(\theta), \theta\in \Theta_1\}\subset I$.
For Bahadur exact slope the following inequality holds:

\begin{equation*}
 c_{T}(\theta)\leq 2K(\theta),
\end{equation*}
where $K(\theta)$ is the Kullback-Leibler information number which measures the statistical distance between the alternative and the
null hypothesis. So, the absolute Bahadur efficiency is defined as

\begin{equation}\label{BE}
 e_{T}(\theta)=\frac{c_{T}(\theta)}{2K(\theta)}.
\end{equation}

In most cases the Bahadur efficiency is not computable for any alternative $\theta$.
However, it is possible to calculate the limit of Bahadur efficiency 
when $\theta$ 
approaches some $\theta_0 \in \Theta_0$. This limit is called the local asymptotic Bahadur efficiency.

Let $G(x;\theta)$ be a family of distributions such that $G(x;0)\in \mathcal{P}$ and $G(x;\theta)\notin \mathcal{P}$ for $\theta\neq 0$.
Then we can reformulate our null hypothesis to be $H_0:\theta=0$. For close alternatives, the local asymptotic Bahadur
efficiency is 

\begin{equation}\label{LBE}
 e_{T}=\lim\limits_{\theta\rightarrow 0}\frac{c_{T}(\theta)}{2K(\theta)}.
\end{equation}

In what follows we shall calculate the local asymptotic Bahadur efficiency for some alternatives and find locally optimal alternatives.
Let $\mathcal{G}=\{G(x;\theta)\}$ be a class of alternatives that satisfy the condition that
it is possible to differentiate along $\theta$ under integral sign in all appearing integrals.
Let $g(x;\theta)$ be a density of a distribution which belongs to $\mathcal{G}$, and let $h(x)=g'_{\theta}(x;0)$. It it easy to see
that $\integ h(x) dx=0$.

\medskip
We now calculate the Bahadur exact slope for our test statistic $T_{n}$. The functions $f(t)$ and $b_{T}(\theta)$
will be determined from the following lemmas.

\begin{lemma}\label{lemmaft}
 Let $t>0$. For statistic $T_n$ the function $f(t)$ is analytic for sufficiently small $t>0$ and it holds
  \begin{equation*}
  f(t)=\frac{54}{5}t^2+o(t^2),\;\;t\rightarrow 0.
 \end{equation*}
\end{lemma}
\textbf{Proof.} Since the kernel $\Upsilon$ is bounded, centered, and non-degenerate, applying the theorem on large
deviations for non-degenerate U-statistics (see \cite{Ponikarov}, Theorem 2.3) we get the statement of the lemma.
$\hfill \Box$

\begin{lemma}\label{lemmaBT}
For a given alternative density $g(x;\theta)$ whose distribution belongs to $\mathcal{G}$ holds
\begin{equation}\label{bT}
  b_T(\theta)=3\theta\integ \upsilon(x)h(x)dx+o(\theta),\;\theta\rightarrow 0.
 \end{equation}
\end{lemma}
\noindent \textbf{Proof:}
From \eqref{tnconv} it follows that
 \begin{eqnarray*}
  \nonumber b_T(\theta)&=&P\left\{\max\left(\frac{X}{Y},\frac{Y}{X}\right)<Z\right\}-\frac{1}{2}
  \\\nonumber  &=&1-P\left\{\frac{X}{Y}>Z,X>Y\right\}-P\left\{\frac{Y}{X}>Z,Y>X\right\}-\frac{1}{2}
  \\\nonumber&=& \frac{1}{2}-2P\left\{\frac{X}{Y}>Z\right\}\\\nonumber&=&
  \frac{1}{2}-2\integ dz\integ dy\int\limits_{yz}^{\infty} g(x;\theta)g(y;\theta)g(z;\theta)dx
  \\&=&\frac{1}{2}-2\integ dz\integ(1-G(yz;\theta)) g(y;\theta)g(z;\theta)dy.
 \end{eqnarray*}
 The first derivative is
 \begin{eqnarray*}
  \nonumber(b_{T})'(\theta)&=&2\integ dz\integ dy \int\limits_{1}^{yz}g'_{\theta}(x;\theta)g(y;\theta)g(z;\theta)dx\\&-&
  4\integ dz \integ(1-G(yz;\theta))g'_{\theta}(y;\theta)g(z;\theta)dy.
 \end{eqnarray*}
 Letting $\theta=0$ we get
\begin{eqnarray*}
\nonumber(b_{T})'(0)&=&2\integ dz\integ dy \int\limits_{1}^{yz}h(x)\alpha y^{-\alpha-1}\alpha z^{-\alpha-1} dx\\\nonumber&-&
  4\integ dz \integ(yz)^{-\alpha}h(y)\alpha z^{-\alpha-1}dy\\\nonumber &=&
  -2\integ dz\integ dy \int\limits_{yz}^{\infty}h(x)\alpha y^{-\alpha-1}\alpha z^{-\alpha-1} dx\\\nonumber&-&
  4\integ dz \integ(yz)^{-\alpha}h(y)\alpha z^{-\alpha-1}dy\\\nonumber &=&
 - 2\integ h(x)dx\int\limits_{1}^x\alpha z^{-\alpha-1}dz\int\limits_{1}^{\frac{x}{z}}\alpha y^{-\alpha-1}dy\\
 &-&
  4\integ h(y)y^{-\alpha}dy\integ \alpha z^{-2\alpha-1}dz\\
  &=&2\integ \alpha x^{-\alpha} h(x)\ln x dx=3\integ \upsilon(x) h(x)dx.
 \end{eqnarray*}
 Since $b_T(0)=0$, using Maclaurin's expansion  for $b_T(\theta)$
 \begin{equation}
  b_T(\theta)=b_T(0)+b'_T(0)\theta+o(\theta),\;\;\theta\rightarrow 0,
 \end{equation}
 we obtain \eqref{bT}.
 $\hfill \Box$

The Kullback-Leibler upper bound $2K(\theta)$ for exact Bahadur slopes for an alternative density $g$ can be determined from the following lemma.

\begin{lemma}\label{lemmaKL}
 For a given density $g(x;\theta)$ let the Kullback-Leibler information
 \begin{equation}\label{KL}
K(\theta)=\inf\limits_{\lambda>0}\int\limits_{1}^\infty\ln\frac{g(x;\theta)}{\lambda x^{-\lambda-1}}g(x;\theta)dx
 \end{equation}
be well-defined. Then when $\theta\rightarrow 0$
\begin{equation}\label{2ka}
 2K(\theta)=\theta^2\bigg(\int\limits_{1}^\infty\frac{ x^{\alpha+1}}{\alpha}h^2(x)dx-
 \bigg(\int\limits_{1}^{\infty}\alpha h(x)\ln x dx\bigg)^2\bigg)+o(\theta^2),
\end{equation}
where $\alpha$ is the shape parameter of the Pareto density $g(x;0)$.
 \end{lemma}
\noindent \textbf{Proof:} The infimum in \eqref{KL} is obtained for $\lambda=
\left(\int\limits_{1}^\infty g(x;\theta)\ln x dx\right)^{-1}$.
Then
\begin{eqnarray*}
\nonumber  K(\theta)&=&\int\limits_{1}^{\infty}g(x;\theta)\ln g(x;\theta)dx+
 \bigg(\frac{1}{\int\limits_{1}^\infty g(x;\theta)\ln x dx}+1\bigg)\int\limits_{1}^\infty g(x;\theta)\ln x dx\\&+&
 \ln \int\limits_{1}^\infty g(x;\theta)\ln x dx\int\limits_{1}^\infty g(x;\theta)dx
\\ &=& \nonumber\integ g(x;\theta)\ln g(x;\theta)dx + 1 + \integ g(x;\theta)\ln x dx + \ln \integ g(x;\theta) \ln x dx.
\end{eqnarray*}
It is obvious from its definition that $K(0)=0$. When we differentiate $K(\theta)$ we get
\begin{equation*}
 K'(\theta)=\int\limits_{1}^\infty g'_{\theta}(x;\theta)\ln g(x;\theta)dx+
 \bigg(1
 +\frac{1}{\int\limits_{1}^\infty g(x;\theta)\ln x dx}\bigg)\int\limits_{1}^\infty g'_{\theta}(x;\theta)\ln x dx.
\end{equation*}
Putting $\theta=0$, we get
\begin{eqnarray*}
 \nonumber K'(0)&=&\int\limits_{1}^\infty h(x)\ln (\alpha x^{-\alpha-1})dx +
 \bigg(1 +\frac{1}{\int\limits_{1}^\infty \alpha x^{-\alpha-1}\ln x dx}\bigg)\int\limits_{1}^\infty h(x)\ln x dx
 \\&=& \ln \alpha \integ h(x)dx - (\alpha +1)\integ h(x)\ln xdx + \integ h(x)\ln x dx
 \\&+& \alpha \integ h(x)\ln x dx = 0.
\end{eqnarray*}
The second derivative of $K(\theta)$ is
\begin{eqnarray*}
 \nonumber K''(\theta)&=&\integ g''_{\theta^2}(x;\theta)\ln g(x;\theta)dx+\integ (g'_{\theta}(x;\theta))^2(g(x;\theta))^{-1}dx\\
 \nonumber &+&\bigg(1+\frac{1}{\integ g(x;\theta)\ln x dx}\bigg)\integ g''_{\theta^2}(x;\theta)\ln x dx\\
 &-&\frac{\bigg(\integ g'_{\theta}(x;\theta)\ln x dx\bigg)^{2}}
 {\bigg(\integ g(x;\theta)\ln x dx \bigg)^{2}}.
\end{eqnarray*}
Putting $\theta=0$, we get
\begin{eqnarray*}
 \nonumber K''(0)&=&\integ g''_{\theta^2}(x;0)\ln (\alpha x^{-\alpha-1})dx+\integ h^2(x)(\alpha x^{-\alpha-1})^{-1}dx\\
 \nonumber &+&\bigg(1+\frac{1}{\integ \alpha x^{-\alpha-1}\ln x dx}\bigg)\integ g''_{\theta^2}(x;0)\ln x dx\\
 &-&\frac{\bigg(\integ h(x)\ln x dx\bigg)^{2}}
 {\bigg(\integ \alpha x^{-\alpha-1}\ln x dx \bigg)^{2}}
 \\&=&\ln \alpha \integ g''_{\theta^2}(x;0) dx-(\alpha+1)\integ g''_{\theta^2}(x;0)\ln xdx+\integ \frac{x^{\alpha+1}}{\alpha}h^2(x)dx
 \\&+&(\alpha+1)\integ g''_{\theta^2}(x;0)\ln x dx -\alpha^2(\integ h(x)\ln x dx )^2 \\
 &=&\int\limits_{1}^\infty\frac{ x^{\alpha+1}}{\alpha}h^2(x)dx-
 \bigg(\int\limits_{1}^{\infty}\alpha h(x)\ln x dx\bigg)^2.
\end{eqnarray*}
From the Maclaurin's expansion of $K(\theta)$
\begin{equation*}
K(\theta)=K(0)+K'(0)\theta+\frac{1}{2}K''(0)\theta^2+o(\theta^2),\;\;\theta\rightarrow 0,
\end{equation*}
 we get \eqref{2ka}.
 $\hfill \Box$

We now present some examples of alternative hypotheses and calculate local asymptotic Bahadur effieciency.

\begin{example}
Let the alternative hypothesis be the log-Weibull distribution with distribution function
\begin{equation*}
 G(x;\theta)=1-e^{-\alpha(\ln x)^{\theta+1}},\;\;x\geq 1, \alpha>0, \theta\in(0,1).
\end{equation*}
The first derivative along $\theta$ of its density at $\theta=0$ is
\begin{equation*}
 h(x)=\frac{\alpha}{x^{\alpha+1}}(-\alpha\ln x\ln \ln x+\ln \ln x+1).
\end{equation*}
Using lemma \ref{lemmaKL}, we get that 
\begin{eqnarray*}
2K(\theta)&=&\theta^2\bigg(\integ \frac{\alpha}{x^{\alpha+1}} ((1-\alpha\ln x)\ln \ln x+1))^2dx
\\&-&
(\integ \frac{\alpha^{2}}{x^{\alpha+1}}((1-\alpha\ln x)\ln \ln x+1) \ln x dx)^{2}\bigg)+o(\theta^{2}),\theta\rightarrow 0,
\end{eqnarray*}
and using lemma \ref{lemmaBT}, we get
\begin{equation*}
 b_T(\theta)=2\theta\integ\frac{\alpha^{2}\ln x}{x^{2\alpha+1}}((1-\alpha\ln x)\ln \ln x+1)dx
 +o(\theta),\theta\rightarrow 0.
\end{equation*}
After calculation of these integrals via expectations of logarithm of gamma distribution and using the property of digamma function
$\psi(x+1)=\psi(x)+\frac{1}{x}$, we obtain
\begin{equation*}
 2K(\theta)=\theta^2\psi'(1)
\end{equation*}
and
\begin{equation*}
 b_{T}(\theta)=\frac{\theta}{4}.
\end{equation*}
From \eqref{slope} and \eqref{LBE}, using lemma \ref{lemmaft}, we get that the local asymptotic Bahadur efficiency is
\begin{equation*}
 e_{T}=\frac{27}{20\psi'(1)}\approx 0.821.
\end{equation*}
 \end{example}

\begin{example}
Let the second alternative hypothesis have the distribution function  given by
\begin{equation*}
 G(x;\theta)=1-e^{-\alpha\ln x-\theta\ln^{\beta}x},\;\; x\geq 1,\;\;\beta>1,\; \theta \in (0,1).
\end{equation*}
Since
\begin{equation*}
 h(x)=\frac{1}{x^{\alpha+1}}(\beta\ln^{\beta-1}x-\alpha\ln^{\beta}x),
\end{equation*}
the Kullback-Leibler bound becomes
\begin{eqnarray*}
 \nonumber 2K(\theta)&=&\theta^2\bigg(\integ\frac{1}{\alpha x^{\alpha+1}}(\beta\ln^{\beta-1}x-\alpha\ln^{\beta}x)^2dx
 \\\nonumber&-&(\integ \ln x \frac{\alpha}{x^{\alpha+1}}(\beta\ln^{\beta-1}x-\alpha\ln^{\beta}x)dx)^2\bigg)+o(\theta^2)
\\&=&\theta^2\frac{\beta^2\Gamma(2\beta-1)-\Gamma^2(\beta+1)}{\alpha^{2\beta}}+o(\theta^2), \theta\rightarrow 0,
\end{eqnarray*}
and
\begin{eqnarray*}
 b_T(\theta)&=&2\theta\integ \frac{\alpha\ln x}{x^{2\alpha+1}}(\beta\ln^{\beta-1}x-\alpha\ln^{\beta}x)dx+o(\theta)\\
 &=&\theta\frac{\Gamma(\beta+1)(\beta-1)}{2^{\beta+1}\alpha^{\beta}}+o(\theta),\theta\rightarrow 0.
\end{eqnarray*}
The local asymptotic Bahadur efficiency is
\begin{equation*}
 e_T=\frac{108}{5}\frac{\Gamma^2(\beta+1)(\beta-1)^2}{2^{2\beta+2}(\beta^2\Gamma(2\beta-1)-\Gamma^2(\beta+1))}.
\end{equation*}
For $\beta=2$ the efficiency is $\frac{27}{80}\approx0.34$. For $\beta=1.5$ the efficiency is
$\frac{27\pi}{160(4-\pi)}\approx 0.62$.
The highest efficiency is reached as $\beta$ approaches 1 and the limit value is $\frac{27}{20\psi'(1)}\approx 0.82$.\\

\end{example}

\subsection{Locally optimal alternative}
Here we study the problem of locally optimal alternatives, the alternatives for which our test statistic
attains the maximal efficiency.
The importance of this problem has been first emphasized in \cite{Bahadur2}. The detailed study was initiated in \cite{Nikitin5}
and developed in \cite{Nikitin1}.
We shall determine some of those alternatives in the following theorem.
\begin{theorem}
Let $\alpha$ be a positive real number and let $g(x;\theta)$ be a density from $\mathcal{G}$ which also satisfies the condition
\begin{equation}
\integ x^{\alpha+1}h^2(x)dx<\infty.
\end{equation}
The alternative densities
\begin{equation*}
 g(x;\theta)=\frac{\alpha}{x^{\alpha+1}}+\theta(C\frac{\alpha\upsilon(x)}{x^{\alpha+1}}+D\frac{\alpha\ln x-1}{x^{\alpha+1}}),\; x\geq 1,\;C>0,\; D\in \mathbb{R},
\end{equation*}
 for small $\theta$ are asymptotically optimal for the test based on $T_n$.
\end{theorem}
\noindent \textbf{Proof:}
Denote
\begin{equation}\label{h0}
 h_0(x)=h(x)-\frac{(\alpha \ln x-1)\alpha^2}{x^{\alpha+1}}\integ h(s)\ln s ds.
\end{equation}
It can be shown  that this function satisfies the following equalities:
\begin{equation}\label{H01}
 \integ h^2_0(x)\alpha^{-1} x^{\alpha+1}dx=\int\limits_{1}^\infty\alpha^{-1} x^{\alpha+1}h^{2}(x)dx-
 \bigg(\int\limits_{1}^{\infty}\alpha h(x)\ln x dx\bigg)^2
\end{equation}
\begin{equation}\label{HO2}
 \integ \upsilon(x)h_0(x)dx=\integ \upsilon(x)h(x)dx.
 \end{equation}
From lemmas \ref{lemmaft} and \ref{lemmaBT}, using   \eqref{sigma}, we get that the
local asymptotic efficiency is
\begin{eqnarray*}
 e_{T}&=&\lim\limits_{\theta\rightarrow 0}\frac{c_T(\theta)}{2K(\theta)}=
 \lim\limits_{\theta\rightarrow 0}\frac{2f(b_T(\theta))}{2K(\theta)}=\lim\limits_{\theta\rightarrow 0}\frac{2\cdot\frac{54}{5}b^2_T(\theta)}{2K(\theta)}
=\lim\limits_{\theta\rightarrow 0}\frac{b_T^2(\theta)}{9\sigma^2 2K(\theta)}\\&=&
\lim\limits_{\theta\rightarrow 0}
 \frac{9\theta^2\bigg(\integ\upsilon(x)h(x)dx\bigg)^2+o(\theta^2)}{9
 \int\limits_{1}^\infty \upsilon^2(x)\alpha x^{-\alpha-1}dx \bigg(\theta^2(\int\limits_{1}^\infty\alpha^{-1} x^{\alpha+1}h^{2}(x)dx-
 \big(\int\limits_{1}^{\infty}\alpha h(x)\ln x dx\big)^2)+o(\theta^{2})\bigg)}\\&=&
 \frac{\bigg(\integ\upsilon(x)h(x)dx\bigg)^2}{
 \int\limits_{1}^\infty \upsilon^2(x)\alpha x^{-\alpha-1}dx \bigg(\int\limits_{1}^\infty\alpha^{-1} x^{\alpha+1}h^{2}(x)dx-
 \big(\int\limits_{1}^{\infty}\alpha h(x)\ln x dx\big)^2\bigg)}
 \\&=&
\frac{\bigg(\integ\upsilon(x)h_0(x)dx\bigg)^2}{
 \int\limits_{1}^\infty \upsilon^2(x)\alpha x^{-\alpha-1}dx \integ h^2_0(x)\alpha^{-1} x^{\alpha+1}dx }.
\end{eqnarray*}
From Cauchy-Schwarz inequality we obtain that $e_{T}=1$ if and only if $h_0(x)=C\upsilon(x)\alpha x^{-\alpha-1}$. Substituting
  this equality in \eqref{h0} we get the expression for $h(x)$. Since $h(x)$ for our alternatives is of such form, we complete the proof.
$\hfill \Box$

\subsection{Critical Values of the Test}
Now we calculate the critical values of this test for small sample sizes.
The statistic $T_{n}$
can be expressed as
\begin{eqnarray*}
 T_{n}&=&\integ (M_{n}(t)-F_{n}(t))dF_{n}(t)= \frac{1}{n}\sum\limits_{j=1}^{n}(M_{n}(x_{j})-F_{n}(x_{j}))\\&=&
 \frac{1}{n}\sum\limits_{j=1}^{n}(\frac{r_{j}-j}{N}-\frac{j}{n})\\&=&
 (2nN)^{-1}(2\sum\limits_{j=1}^{n}r_{j}-(n+1)(n+N)),
\end{eqnarray*}
where $N=\binom{n}{2}$, and $r_{j}$ is the rank of ${x_{j}}$ in the
pooled sample of $x_{k}$, $1\leq k \leq n$ and $\max\{\frac{x_{j}}{x_{k}},\frac{x_{k}}{x_{j}}\}$, $1\leq j<k \leq n$.

Since we don't have exact distribution for small values of $n$, the critical values of the test can be calculated using Monte-Carlo methods.
The one-tailed critical values, based on 10000 repetitions, are given in table \ref{fig: cvTn}.
\begin{table}[ht]
\centering

\begin{tabular}{r|rrr}
 & \multicolumn{3}{c}{level of significance} \\

  $n $& 0.1 & 0.05 & 0.01 \\
  \hline
10 & 0.09 & 0.13 & 0.22 \\
  20 & 0.06 & 0.09 & 0.15 \\
  30 & 0.05 & 0.07 & 0.11 \\
  40 & 0.04 & 0.06 & 0.09 \\
  50 & 0.04 & 0.06 & 0.09 \\
  100 & 0.03 & 0.04 & 0.05 \\
   \hline
\end{tabular}
\caption{Critical values for the statistic $T_n$}
\label{fig: cvTn}
\end{table}

\newpage
\section{Statistics $V_n$}

In this section we examine the asymptotic properties of Kolmogorov-Smirnov type statistic $V_n$ under null hypothesis.
For fixed $t\in[1,\infty)$ the expression $M_n(t)-F_n(t)$ is an 
$U$-statistic with kernel
\begin{equation*}
  \varPsi(X,Y;t)=I\left\{\max\left(\frac{X}{Y},\frac{Y}{X}\right)\leq t\right\}-
  \frac{1}{2}I\left\{X\leq t\right\}-\frac{1}{2}I\left\{Y\leq t\right\}.
\end{equation*}
Let $\psi(X;t)$ be the projection of $\varPsi(X,Y;t)$ on $X$. Then
\begin{eqnarray*}
\psi(s;t)&=& E(\varPsi(X,Y;t)|X=s)\\&=&P\left\{\max\left(\frac{s}{Y},\frac{Y}{s}\right)\leq t\right\}
-\frac{1}{2}I\{s\leq t\}-\frac{1}{2}P\left\{Y\leq t\right\}
 \\&=&\frac{t^{\alpha}}{s^{\alpha}}-\frac{1}{s^{\alpha} t^{\alpha}}-
 \frac{1}{2}+\frac{1}{2t^{\alpha}}+I\{s\leq t \}\left(\frac{1}{2}-\frac{t^{\alpha}}{s^{\alpha}}\right).
\end{eqnarray*}
It is easy to show that expected value of  $\psi(X;t)$ is zero.
Its variance for fixed $t$ is
\begin{equation*}
\sigma^2(t)=Var(\psi(X;t))=\frac{1}{12t^{3\alpha}}(t^{2\alpha}+t^{\alpha}-2).
\end{equation*}
The function $\sigma^2(t)$ reaches its maximum for $t_0=(\sqrt{7}-1)^\frac{1}{\alpha}$,
and that maximum is equal to $\frac{7\sqrt{7}+10}{648}\approx 0.044$. Hence our family of kernels  $\varPsi(X,Y;t)$
using argumentation from \cite{Nikitin} is not degenerate.
It can be shown using \cite{silverman} that U-empirical random process
$\rho(t)=\sqrt{n}(M_n(t)-F_n(t)),\;t\geq 1$, converges
in distribution to some Gaussian process. It is not easy to calculate the covariance of this process, and the asymptotic distribution
of statistic $V_n$ is unknown.

\subsection{Bahadur Efficiency}
 Now we shall calculate Bahadur efficiency in the analogous way as in previous section. Here, the function $f(t)$
 from \eqref{slope} for statistic $V_n$ is determined  in  the following theorem.
\begin{theorem}\label{theo}
Let $t>0$. Then $f(t)$ is analytic for sufficiently small $t>0$ and it holds
 \begin{equation*}
  f(t)=\frac{t^2}{8\sigma^2(t_0)}+o(t^2)\sim 2.84 t^2,\;\;t\rightarrow 0.
 \end{equation*}
\end{theorem}
The proof of this theorem can be found in \cite{Nikitin}.

\medskip
\noindent In the following lemma we determine $b_V(\theta)$, the limit in probability of $V_n$.
\begin{lemma}\label{lemmabV}
For a given alternative density $g(x;\theta)$ whose distribution belongs to $\mathcal{G}$ holds
\begin{equation}\label{bV}
  b_V(\theta)=2\theta\sup\limits_{t\geq 1}\left|\integ \psi(x;t)h(x)dx\right|+o(\theta),\;\theta\rightarrow 0.
 \end{equation}
\end{lemma}
\textbf{Proof.}
Using Glivenko-Cantelli theorem for $U$-empirical distribution functions (see \cite{jansen}), we get

\begin{eqnarray}\label{b0}
\nonumber b_V(\theta)&=&\sup\limits_{t\geq 1}\left|P\left\{\max\left(\frac{X}{Y},\frac{Y}{X}\right)\leq t\right\}-G(t;\theta)\right|\\
 &=&\sup\limits_{t\geq 1}\left|2\integ dx\int\limits_{x}^{xt}g(y;\theta)g(x;\theta)dy-\int\limits_{1}^tg(x;\theta)dx\right|.
\end{eqnarray}
Let us denote
\begin{equation*}
 a_V(\theta)=2\integ dx\int\limits_{x}^{xt}g(y;\theta)g(x;\theta)dy-\int\limits_{1}^tg(x;\theta)dx.
\end{equation*}
It is easy to show that $a_V(0)=0$.

\noindent The first derivative of $a_V(\theta)$ is
\begin{eqnarray*}
a_V'(\theta)&=&2\integ dx\int\limits_{x}^{xt}g'_{\theta}(y;\theta)g(x;\theta)dy+
2\integ dx\int\limits_{x}^{xt}g(y;\theta)g'_{\theta}(x;\theta)dy\\&-&\int\limits_{1}^tg'_{\theta}(x;\theta)dx.
\end{eqnarray*}

\noindent The first derivative at $\theta=0$ is

\begin{eqnarray*}
 \nonumber a_V'(0) &=& 2\integ dx\int\limits_{x}^{xt}h(y)\frac{\alpha}{x^{\alpha+1}}dy+2\integ dx\int\limits_{x}^{xt}\frac{\alpha}{y^{\alpha+1}}h(x)dy-
\int\limits_{1}^th(x)dx
\\\nonumber &=& \int\limits_{1}^th(x)dx+2t^\alpha\int\limits_{t}^{\infty}h(x)x^{-\alpha}dx-2t^{-\alpha}\integ h(x)x^{-\alpha} dx
\\\nonumber &=&\integ\left[I\left\{x\leq t\right\}+2t^{\alpha}x^{-\alpha}\left(1-I\left\{x\leq t\right\}\right)-2t^{-\alpha}x^{-\alpha}\right]h(x)dx.
\end{eqnarray*}

\noindent Applying Maclaurin's expansion to $a_V(\theta)$
\begin{equation*}
  a_V(\theta)=a_V(0)+a'_V(0)\theta+o(\theta),\;\;\theta\rightarrow 0,
 \end{equation*}
 and inserting this expression in \eqref{b0}
we obtain \eqref{bV}.
$\hfill \Box$

\begin{example}
 Let the alternative hypothesis be the mixture of two Pareto distributions with the following distribution function
\begin{equation}\label{mix}
 G(x;\theta)=(1-\theta)(1-x^{-\alpha})+\theta(1-x^{-\beta}), \;\;x\geq 1, \beta>\alpha>0, \theta\in(0,1).
 \end{equation}
 The first derivative along $\theta$ of its density at $\theta=0$ is
\begin{equation*}
 h(x)=-\frac{\alpha}{x^{\alpha+1}}+\frac{\beta}{x^{\beta+1}}.
\end{equation*}
Using lemma \ref{lemmaKL}, we get that
\begin{eqnarray*}
\nonumber2K(\theta)&=&\theta^2\bigg(\integ \frac{x^{\alpha+1}}{\alpha}(-\frac{\alpha}{x^{\alpha+1}}
+\frac{\beta}{x^{\beta+1}})^2dx\\&-&
(\integ \alpha
(-\frac{\alpha}{x^{\alpha+1}}+\frac{\beta}{x^{\beta+1}})\ln x dx)^2\bigg)+o(\theta^{2})
\\\nonumber&=&\theta^2\frac{(\beta-\alpha)^4}{\alpha\beta^2(2\beta-\alpha)}
+o(\theta^{2}), \theta\rightarrow 0,
\end{eqnarray*}
and using lemma \ref{lemmabV}, we get
\begin{eqnarray*}
 \nonumber b_V(\theta)&=&2\theta \sup\limits_{t\geq 1}\left|\integ\left(I\left\{x\leq t\right\}(\frac{1}{2}-\frac{t^\alpha}{x^\alpha})+x^{-\alpha}
 (t^{\alpha}-t^{-\alpha})\right)(\frac{-\alpha}{x^{\alpha+1}}+\frac{\beta}{x^{\beta+1}})dx\right|
 \\ &=& 2\theta \sup\limits_{t\geq 1}\left|(t^{-\alpha}-t^{-\beta})(\frac{\beta}{\alpha+\beta}-\frac{1}{2})\right|
 +o(\theta)
 \\&=& \theta \frac{(\beta-\alpha)^2}{\alpha(\alpha+\beta)(\frac{\beta}{\alpha})^{\frac{\beta}{\beta-\alpha}}}
 +o(\theta), \;\theta\rightarrow 0.
 \end{eqnarray*}
 Using \eqref{slope}, \eqref{LBE} and theorem \ref{theo} we calculate the local asymptotic Bahadur efficiency. We get
 \begin{equation}
  e_V=\lim\limits_{\theta\rightarrow 0}\frac{c_V(\theta)}{2K(\theta)}=5.68\frac{\alpha(2\beta-\alpha)}{(\alpha+\beta)^2(\frac{\beta}{\alpha})^{\frac{2\alpha}{\beta-\alpha}}}.
 \end{equation}
 This expression reaches its maximum for $\beta=4.646\;\alpha$ and then
$e_V\approx 0.636$.
\end{example}

\subsection{Locally optimal alternative}
As in the previous section,  we shall determine some of locally optimal alternatives in the following theorem.
\begin{theorem}
Let $\alpha$ be a positive real number and let $g(x;\theta)$ be a density from $\mathcal{G}$ which also satisfies the condition
\begin{equation}
\integ x^{\alpha+1}h^2(x)dx<\infty.
\end{equation}
The alternative densities
\begin{equation*}
 g(x;\theta)=\frac{\alpha}{x^{\alpha+1}}+\theta(C\frac{\alpha\psi(x;t_0)}{x^{\alpha+1}}+D\frac{\alpha\ln x-1}{x^{\alpha+1}}),\; x\geq 1,\;C>0,\; D\in \mathbb{R},
\end{equation*}
 where $t_0=(\sqrt{7}-1)^{\frac{1}{\alpha}}$, for small $\theta$ are asymptotically optimal for the test based on $V_n$.
\end{theorem}
\textbf{Proof:}
Let $h_0$ be the function defined in \eqref{h0}.
It can be shown  that this function besides \eqref{H01}, also satisfies the following equality:
\begin{equation}\label{HO2V}
 \integ \psi(x;t)h_0(x)=\integ \psi(x;t)h(x).
 \end{equation}
From theorem \ref{theo} and lemma \ref{lemmabV}, we get that the
local asymptotic efficiency is
\begin{eqnarray*}
 e_V&=&\lim\limits_{\theta\rightarrow 0}\frac{c_V(\theta)}{2K(\theta)}=
 \lim\limits_{\theta\rightarrow 0}\frac{2f(b_V(\theta))}{2K(\theta)}
=\lim\limits_{\theta\rightarrow 0}\frac{b_V^2(\theta)}{4\sigma^2(t_0) 2K(\theta)}\\&=&
 \lim\limits_{\theta\rightarrow 0}
 \frac{4\theta^2\sup\limits_{t\geq 1}\bigg(\integ\psi(x;t)h(x)dx\bigg)^2+o(\theta^2)}{4
 \sup\limits_{t\geq 1}\int\limits_{1}^\infty \psi^2(x;t)\alpha x^{-\alpha-1}dx \bigg(\theta^2(\int\limits_{1}^\infty\alpha^{-1} x^{\alpha+1}h^{2}(x)dx-
 \big(\int\limits_{1}^{\infty}\alpha h(x)\ln x dx\big)^2)+o(\theta^{2})\bigg)}\\&=&
 \frac{\sup\limits_{t\geq 1}\bigg(\integ\psi(x;t)h(x)dx\bigg)^2}{\sup\limits_{t\geq 1}
 \int\limits_{1}^\infty \psi^2(x;t)\alpha x^{-\alpha-1}dx \bigg(\int\limits_{1}^\infty\alpha^{-1} x^{\alpha+1}h^{2}(x)dx-
 \big(\int\limits_{1}^{\infty}\alpha h(x)\ln x dx\big)^2\bigg)}
 \\&=&
\frac{\sup\limits_{t\geq 1}\bigg(\integ\psi(x;t)h_0(x)dx\bigg)^2}{
\sup\limits_{t\geq 1} \int\limits_{1}^\infty \psi^2(x;t)\alpha x^{-\alpha-1}dx \integ h^2_0(x)\alpha^{-1} x^{\alpha+1}dx }.
\end{eqnarray*}
From Cauchy-Schwarz inequality we obtain that $e_V=1$ if $h_0(x)=C\psi(x;t_0)\alpha x^{-\alpha-1}$. Substituting
  this equality in \eqref{h0} we get the expression for $h(x)$. Since $h(x)$ for our alternatives is of such form, we complete the proof.
$\hfill \Box$

\subsection{Critical Values}
Since we do not know the distribution of $V_n$,
the critical values of the test can be calculated using Monte-Carlo methods.
The two-tailed critical values, based on 10000 repetitions, are given in table \ref{fig: cvVn}.
\begin{table}[ht]
\centering
\begin{tabular}{r|rrr}
& \multicolumn{3}{c}{level of significance} \\
$n$ & 0.1 & 0.05 & 0.01 \\
\hline
10 & 0.36 & 0.41 & 0.56 \\
20 & 0.24 & 0.28 & 0.36 \\
30 & 0.19 & 0.22 & 0.28 \\
40 & 0.16 & 0.19 & 0.24 \\
50 & 0.15 & 0.16 & 0.21 \\
100 & 0.11 & 0.12 & 0.15 \\
\hline
\end{tabular}
\caption{Critical values for the statistic $V_n$}
\label{fig: cvVn}
\end{table}

\section{Power Comparison}
In this section we compare the powers of our tests with the powers of two tests most commonly used for these types of hypotheses,
namely
Kolmogorov-Smirnov ($D_{n}$) and Cramer-von Mises ($\omega^{2}_{n}$) tests. The comparison is done for sample sizes of 20 and 50
for level of significance 0.05.
We propose five distributions that are usually considered alternatives for Pareto distribution:
\begin{itemize}
 \item log-normal with $m=0,\sigma=1$
 \item half-normal with $\sigma=1$
 \item Weibull with  $\alpha=2$
 \item gamma with  $\alpha=2$, $\beta=1$
 \item log-gamma with $\alpha=2,\beta=1$.
\end{itemize}
For Kolmogorov Smirnov and Cramer von Mises tests, since they are not applicable to composite hypothesis,
we first estimated shape parameter $\alpha$ with its MLE  $\hat{\alpha}=n(\sum\limits_{k=1}^n \ln X_k)^{-1}$ and we
calculate critical values of corrected test using Monte Carlo procedure. We then  calculate powers for sample sizes $n=20$ and $n=50$, using Monte Carlo method.
The powers are given in table \ref{fig: moc1}.

We can notice that for all given alternatives and for both sample sizes our tests based on the statistics $T_n$ and $V_n$ 
 have greater powers than modified Kolmogorov-Smirnov test. In comparison with modified Cramer-von Mises test
 our tests have similar powers, and in most cases at least one of our tests performs better.
\begin{center}
\begin{table}[ht]

\centering
\begin{tabular}{|c|c|c|c|c|c|}\hline
 $n$ & alternative & $T_{n}$ & $V_{n}$ & $D_{n}$ & $\omega^{2}_{n}$ \\ \hline
 \multirow{5}{5mm}{20} & log-normal & 0.6263 & 0.6713 & 0.5585 & 0.6432\\
 & halfnormal & 0.6254 & 0.6718 & 0.5327 & 0.6489\\
 & Weibull & 0.9984 & 0.9988 & 0.9893 & 0.9990 \\
 & gamma & 0.9937 & 0.9940 & 0.9642 & 0.9919\\
 & log-gamma & 0.4654 & 0.5282 & 0.4096 & 0.4643 \\\hline
 \multirow{5}{5mm}{50} & log-normal & 0.9877 & 0.9758 & 0.9520 & 0.9841 \\
 & halfnormal & 0.9697 &
 0.9691 & 0.9268 & 0.9763 \\
 & Weibull & 1 & 1 & 1 & 1\\
 & gamma & 1 & 1 & 1 & 1\\
 & log-gamma & 0.9158 & 0.8955 & 0.8241 & 0.9015\\\hline
\end{tabular}
\caption{Power of the test for different alternatives}
\label{fig: moc1}
\end{table}
\end{center}

\section{Conclusion}

In this paper we gave a new characterization of Pareto distribution and proposed two goodness of fit tests based on it.

The main advantage of our tests is the fact that they are free of parameter $\alpha$ which enables us to test a composite null
hypothesis.

The Bahadur efficiency for some alternatives has been calculated and the obtained efficiencies are reasonably high. 
For both tests we determined locally optimal class of alternatives. We compared these tests with some commonly used goodness
of fit tests and we noticed that in most cases at least one of our tests performs slightly better. 

\section*{Acknowledgement}

We would like to thank the editor and the anonymous referees for their very useful comments and suggestions that greatly 
improved our paper.

\end{document}